\def\overrideversion{\defaultversion}

\newif\ifdcg
\dcgtrue

\newif\ifsocg
\socgfalse

\newif\ifexperiment
\experimentfalse

\ifdefined\overrideversion
    \socgfalse
    \dcgfalse
\else
    \socgfalse
    \dcgtrue
\fi

\ifsocg
\documentclass[a4paper,UKenglish,cleveref,nosubfigcap,numberwithinsect]{socg-lipics-v2019}
\author{Brett Leroux}{Department of Mathematics, University of California, Davis,  USA}{beleroux@ucdavis.edu}{}{(Optional) author-specific funding acknowledgements}
\author{Luis Rademacher}{Department of Mathematics, University of California, Davis,  USA \and \url{https://www.math.ucdavis.edu/~lrademac/} }{lrademac@ucdavis.edu}{}{(Optional) author-specific funding acknowledgements}
\authorrunning{B. Leroux and L. Rademacher}
\Copyright{Brett Leroux and Luis Rademacher}
\category{}
\relatedversion{}
\supplement{}

\newtheorem{problem}[theorem]{Problem}
\newtheorem{conjecture}[theorem]{Conjecture}

\ifexperiment
\usepackage[location=appendix,appendixsectionname={Omitted Proofs}]{moveproofs}
\usepackage{environ}
\else
\usepackage[location=appendix,appendixsectionname={Omitted Proofs}]{moveproofs}
\renewcommand{\appendixproofsection}[1]{\section{#1}\label{sec:omittedproofs}}
\DeclareDocumentCommand{\appendixproof}{s m}{
\ifmoveproofstoappendix
    \IfBooleanTF{#1}{
        \MPGetProof{#2}
    }{
        \begin{proof} (of \cref{#2})
        \MPGetProof{#2}
        \end{proof}
    }
\else
\fi
}

\else

\ifdcg
\RequirePackage{fix-cm}
\documentclass[smallextended,numbook,envcountsame,envcountreset,nospthms]{svjour3}
\usepackage{graphicx}
\usepackage{amsfonts,amsmath,amssymb,amscd,enumitem,color}
\usepackage{mathtools}
\usepackage[hidelinks]{hyperref}
\makeatletter
\let\cl@chapter\undefined
\makeatletter
\usepackage[english]{babel}
\usepackage[utf8]{inputenc}

\usepackage{amsthm}
\usepackage[capitalize]{cleveref}
\newtheorem{theorem}{Theorem}
\newtheorem{proposition}[theorem]{Proposition}
\newtheorem{definition}[theorem]{Definition}
\newtheorem{lemma}[theorem]{Lemma}
\newtheorem{problem}[theorem]{Problem}
\newtheorem{conjecture}[theorem]{Conjecture}
\newtheorem{example}[theorem]{Example}
\numberwithin{theorem}{section}

\author{Brett Leroux \and Luis Rademacher}

\institute{Brett Leroux \at
              University of California, Davis, CA, USA \\
              \email{beleroux@ucdavis.edu}           
           \and
           Luis Rademacher \at
              University of California, Davis, CA, USA \\
              \email{lrademac@ucdavis.edu}
}

\date{Received: date / Accepted: date}

\newcommand{\makeproof}[3]{\begin{proof}{#3}\end{proof}}

\else

\documentclass[11pt]{article}
\usepackage{graphicx}
\usepackage{amsfonts,amsmath,amssymb,amscd,enumitem,color}
\usepackage{mathtools}
\usepackage{accents}
\usepackage{amsthm}
\usepackage[hidelinks]{hyperref}
\usepackage[capitalize]{cleveref}
\usepackage[english]{babel}
\usepackage[utf8]{inputenc}

\newtheorem{theorem}{Theorem}
\newtheorem{proposition}[theorem]{Proposition}
\newtheorem{definition}[theorem]{Definition}
\newtheorem{lemma}[theorem]{Lemma}
\newtheorem{problem}[theorem]{Problem}
\newtheorem{conjecture}[theorem]{Conjecture}
\newtheorem{example}[theorem]{Example}
	
\numberwithin{theorem}{section}

\usepackage{fullpage}

\newenvironment{acknowledgements}{\paragraph{Acknowledgements}}{}

\author{Brett Leroux, Luis Rademacher}

\def\keywords#1{\par\addvspace\medskipamount{\rightskip=0pt plus1cm
\def\and{\ifhmode\unskip\nobreak\fi\ $\cdot$
}\noindent{Keywords:}\enspace\ignorespaces#1\par}}

\newcommand{\makeproof}[3]{\begin{proof}{#3}\end{proof}}
\fi
\fi

\title{Algebraic \texorpdfstring{\(k\)}{k}-sets and generally neighborly embeddings}

\usepackage{dirtytalk}

\def\final{1}  

\ifnum\final=0  
\newcommand{\lnote}[1]{[{\small Luis: \bf #1}]}
\newcommand{\bnote}[1]{[{\small Brett: \bf #1}]}
\newcommand{\anonnote}[1]{[{\small anon: \bf #1}]}
\newcommand{\sidecomment}[1]{\marginpar{\tiny #1}}
\newcommand{\details}[1]{{\color{blue}\ [[#1]] }}
\else 
\newcommand{\lnote}[1]{}
\newcommand{\bnote}[1]{}
\newcommand{\anonnote}[1]{}
\newcommand{\sidecomment}[1]{}
\newcommand{\details}[1]{}
\fi  

\newcommand{\Rl}{\operatorname{\mathbb{R}}}
\newcommand{\Nl}{\operatorname{\mathbb{N}}}
\newcommand{\RR}{\mathbb{R}}
\newcommand{\relint}{\operatorname{relint}}

\newcommand{\Cx}{\operatorname{\mathbb{C}}}
\newcommand{\F}{\operatorname{\mathcal{F}}}

\newcommand{\Tor}{\operatorname{\Tor}}

\newcommand{\aff}{\operatorname{aff}}
\newcommand{\conv}{\operatorname{conv}}
\newcommand{\interior}{\operatorname{int}}
\newcommand{\suchthat}{\mathrel{:}}
\newcommand{\bd}{\operatorname{bd}}

\makeatletter
\newcommand\fake@math{}
\def\fake@math#1\){[math]}
\makeatother

\begin{document}
\maketitle

\ifexperiment
%

\let\oldlabel\label
\renewcommand{\label}[1]{%
\gdef\currentlabel{#1}%
\oldlabel{#1}}
\RenewEnviron{proof}{\globaldefs=1\makeproof{a}{\currentlabel}{\BODY}}
\fi

\begin{abstract}
Given a set $S$ of $n$ points in $\Rl^d$, a $k$-set is a subset of $k$ points of $S$ that can be strictly separated by a hyperplane from the remaining $n-k$ points. Similarly, one may consider $k$-facets, which are hyperplanes that pass through $d$ points of $S$ and have $k$ points on one side. 
A notorious open problem is to determine the asymptotics of the maximum number of $k$-sets.
In this paper we study a variation on the $k$-set/$k$-facet problem with hyperplanes replaced by algebraic surfaces.
In stark contrast to the original $k$-set/$k$-facet problem, there are some natural families of algebraic curves for which the number of $k$-facets can be counted exactly. 
For example, we show that the number of halving conic sections for any set of $2n+5$ points in general position in the plane is $2\binom{n+2}{2}^2$. 
To understand the limits of our argument we study a class of maps we call \emph{generally neighborly embeddings}, which map generic point sets into neighborly position. 
Additionally, we give a simple argument which improves the best known bound on the number of $k$-sets/$k$-facets for point sets in convex position. 

\end{abstract}

\section{Introduction}

The $k$-set problem aims to understand the following: given a set of $n$ points in $\RR^d$, in how many ways can it be strictly separated into $k$ and $n-k$ points by an affine hyperplane? 
More specifically, one wants to understand the asymptotics of the maximum number of such partitions as a function of $n$ and $d$.
The question of determining the maximum number of $k$-sets for point sets in $\Rl^2$ was first raised by A.~Simmons in unpublished work. Straus, also in unpublished work, gave a construction showing an $\Omega(n\log n)$ lower bound. Lov\'{a}sz \cite{lovasz1971number} published the first paper on $k$-sets, establishing an $O(n^{3/2})$ upper bound. See also the paper \cite{MR0363986} of Erd\H{o}s, Lov\'{a}sz, Simmons,  and Straus.

The main challenge is that even for the basic $k$-set problem on the plane, the asymptotics of the maximum number of $k$-sets is not well understood despite decades of effort.
The best known asymptotic bounds are $n e^{\Omega(\sqrt{\log k})}$ \cite{DBLP:conf/compgeom/Toth00,Tothlower,MR2405686} and $O(nk^{1/3})$ \cite{Dey1997ImprovedBO,DBLP:journals/dcg/Dey98}.

It is natural to ask similar questions for families of surfaces different from all hyperplanes.
These sorts of questions have been studied in \cite{1676031,Ardila04,DBLP:conf/compgeom/BuzagloHP08,DBLP:conf/compgeom/ChevallierFSS14,Chevallier}.
Ardila's paper \cite{Ardila04} shows that for any set of $2n + 1$ points in general position in the plane and any $0 \le k \le 2n-2$, the number of circles that go through 3 points and have $k$ points on one side is exactly $2(k+1)(2n-k-1)$.
We call this phenomenon \emph{exact counting}: when for a family of curves (or surfaces), there exists an integer $d$ such that given any generic set of $n$ points and any $k$, the number of curves which pass through $d$ points and have $k$ points on one side depends only on $n$ and $k$ and not on the points.
A result essentially equivalent to Ardila's was proven earlier in \cite{1676031} by counting vertices of certain Voronoi diagrams. Chevallier et al.\ extended the result to convex pseudo-circles in \cite{Chevallier}.

Borrowing the language of set theory/computational geometry/learning theory, one can think of the $k$-set problem as being formulated over a \emph{set system} (also known as a \emph{hypergraph}, \emph{hypothesis class} or \emph{range space}), namely a universe and a family of subsets of the universe. 
In the $k$-set problem the universe is $\RR^d$ and the family of subsets is all halfspaces.
This paper takes a step towards the understanding of the $k$-set problem for general set systems.
We focus on set systems induced by maps in the following way: given a map $\varphi:\RR^d \to \RR^p$, the set system induced by $\varphi$ has universe $\RR^d$ and family of subsets $\{ \varphi^{-1}(H) \suchthat \text{$H$ is a closed halfspace in $\RR^p$}\}$. 
Moreover, most of our results involve maps $\varphi$ with components that are polynomials, so that the separating surfaces in the resulting set system are algebraic surfaces.
One of our main examples is the Veronese map (\cref{def:Veronese}) which induces separators that are algebraic surfaces of degree at most $m$. 
The Veronese map is also known as the feature map of the polynomial kernel in machine learning \cite{Shawe-Taylor:2004:KMP:975545}.




\subsection*{Our contributions:}
\begin{itemize}
\item
\textbf{Exact count.}
We show that the exact count phenomenon of \cite{1676031,Ardila04} (for halving circles) holds for other natural set systems: 
conic sections (\cref{thm:coniccount}) and homogeneous polynomials of fixed even degree on the plane (\cref{thm:homogeneouscount}). 
We prove this by establishing a remarkable property of the corresponding maps: 
generic point sets are mapped to point sets that form the vertices of a neighborly polytope (\cref{prop:ConicNeighborly,prop:HomNeighborly}, see \cref{sec:neighborly} for background).
This is then combined with the known fact that the number of $k$-facets of a neighborly point set is given by a formula that depends only on the dimension, on $k$, and on the number of points \cite{Clarkson87,DBLP:journals/dcg/ClarksonS89}, {\cite[Proposition 4.1]{WagnerSurvey}}.


\item 
\textbf{Limits of the neighborliness argument.}
We study the limits of the neighborliness argument above that provides exact counting. 
We show that, for maps whose image is a variety, the argument only works for points on the plane.
We proceed as follows:
For the argument to work, one needs the map $\varphi: \RR^d \to \RR^p$ to map a generic set of points into a $k$-neighborly set of points for certain $k$. When $\varphi$ is an embedding, we call the image $M:=\varphi(\Rl^d)$ a \emph{generally $k$-neighborly $d$-manifold} (\cref{def:generallyneighborlymanifold}).
We study the minimal dimension $p$ so that $M$ is a generally $k$-neighborly $d$-manifold and show that $p \leq 2k+d-1$ (\cref{thm:embedding}). For the same question with manifolds replaced by algebraic varieties, we show that $p = 2k+d-1$ (\cref{prop:lowerbound}).
This line of work relates to a problem of M.~Perles on $k$-neighborly embeddings (see \cref{subsec:perles}). 



\item
\textbf{Convex position bound.} 
We show an improved upper bound on the number of $k$-sets/$k$-facets for points in convex position (\cref{prop:convex}). 
While our argument is simple, we are not aware of any known bounds for the convex case better than the general case in dimension higher than 3 (the convex case is well understood in two and three dimensions).
\item
\textbf{Degree of neighborliness.}
We study the degree of neighborliness of point sets mapped by a $\varphi$ with components ``all monomials of degree at most $m$'' or ``all monomials of degree exactly $m$'' (\cref{prop:veroneseneighborly,prop:homogeneousneighborly}). 
In particular, for even $m$, point sets are mapped into point sets in convex position and the convex position bound gives an improved bound on the number of $k$-facets.

\item
\textbf{Weakly neighborly point sets.} 
We leverage \emph{weakly $k$-neighborly} point sets (\cref{def:weakly}), a notion that is better behaved for our purposes than generally $k$-neighborly and Perles' $k$-neighborly maps (\cref{prop:compactness}). 
In particular, we study weakly $k$-neighborly algebraic varieties, and resolve the question of the minimal $p$ so that $\Rl^p$ contains a weakly $k$-neighborly $d$-dimensional algebraic variety. We show that the minimal dimension is $2k+d-1$ (\cref{thm:lowerboundweak}).
\end{itemize}

\paragraph{Outline of the paper.} \cref{sec:preliminaries} introduces the $k$-facet problem and the generalization to set systems which are induced by maps. In \cref{sec:count} we count $k$-facets for set systems induced by maps. This amounts to studying $k$-facets of point sets of the form $\varphi(S)$ where $\varphi: \Rl^d \to \Rl^p$ is some map and $S \subset \Rl^d$ is a finite set of points.  \cref{sec:limits} concerns the limits of the neighborliness argument mentioned above. We conclude with some open questions in \cref{sec:conclusion}.
 
\ifsocg
Omitted proofs can be found in \cref{sec:omittedproofs}.
\fi




\section{Preliminaries}\label{sec:preliminaries}

Let $S$ be a set of points in $\Rl^d$. A \emph{$k$-set} of $S$ is a subset $A \subset S$ of size $k$ that can be strictly separated from $S \setminus A$ by a hyperplane. 
\begin{definition}
Let $a_k(S)$ be the number of $k$-sets of the set $S \subset \Rl^d$ and $a_k^{(d)}(n)$ be the maximum number of $k$-sets that a set of $n$ points in $\Rl^d$ can have. 
\end{definition}

When studying the $k$-set problem, one usually only considers point sets which are in general linear position. 
\begin{definition}\label{def:genlinear}
A set of at least $d+1$ points in $\Rl^d$ is in \emph{general linear position} if no $d+1$ (and thus, fewer) points are affinely dependent. 
\end{definition}
This reduction is justified by the observation that the maximum number of $k$-sets is attained by a set of points in general linear position (see for example \cite{WagnerSurvey}). 
For point sets in general linear position, one can study the closely related concept of \emph{$k$-facets}.

\begin{definition}
Let $S$ be a finite set of points in general linear position in $\Rl^d$ and let $\Delta$ be a subset of $d$ points from $S$. The subset $\Delta$ along with some orientation of the hyperplane $\aff \Delta$ is a \emph{$k$-facet} of $S$ if the open halfspace on the positive side of $\aff \Delta$ contains exactly $k$ points from $S$. 
\end{definition}

\begin{definition}
Let $e_k(S)$ be the number of $k$-facets of the set $S \subset \Rl^d$ and $e_k^{(d)}(n)$ be the maximum number of $k$-facets that a set of $n$ points in general linear position in $\Rl^d$ can have. 
\end{definition}

It seems unlikely that one would be able to determine $e_k^{(d)}(n)$ or $a_k^{(d)}(n)$ precisely, so instead efforts have focused on finding the asymptotic behavior of these functions. 
If one is only concerned with the asymptotics, then it suffices to study either $k$-sets or $k$-facets since for fixed $d$ and $n \to \infty$, $a_k^{(d)}(n)$ and $e_k^{(d)}(n)$ have the same asymptotic behavior \cite{WagnerSurvey}. 
In this paper we focus on $k$-facets since our results which count $k$-facets exactly cannot immediately be adapted to the setting of $k$-sets.

\subsection{Generic properties and general position}\label{sec:generic}
After defining $k$-sets and $k$-facets for set systems other than halfspaces, we will need to use various notions of general position different from general linear position. 

A \emph{generic property} (of point sets) is one that holds for all but a relatively small number of atypical point sets.
The point sets which satisfy a generic property are said to be \emph{in general position}. 
We use the terms ``generic point set'' and ``point set in general position'' interchangeably. 


In algebraic geometry a generic property is one that holds for a dense and open set. In other fields, a generic property holds almost everywhere.
For concreteness in some of our statements we set ``generic property'' to mean a property that holds in an open and dense set, but this choice is not always crucial. 


\subsection{Set systems}\label{sec:setsystems}
A set system is a pair $(X,\F)$ where $X$ is a ground set (or universe) and $\F$ is a collection of subsets of $X$. 

We will restrict our attention to set systems which are \emph{induced by maps}. Suppose we have a map $\varphi: \Rl^d \to \Rl^p$, that is, a map of $\Rl^d$ into some (usually higher dimensional) space. Any such map induces a set system on the ground set $\Rl^d$ in the following way. 
Let $\F_\varphi$ consist of all regions $R \subset \Rl^d$ of the form $\varphi^{-1}(H )$ where $H$ is a closed halfspace in $\Rl^p$. 
We say that $R$ is induced by the halfspace $H$ and we say that the set system $(\Rl^d,\F_\varphi)$ is induced by $\varphi$. Many interesting set systems are induced by maps.

\subsection{Set systems induced by Veronese-type maps}
\begin{definition}
A \emph{polynomial map} is a map $\Rl^d \to \Rl^p$ defined by 
\[
x \mapsto (f_1(x), \dotsc, f_p(x))
\] 
where $f_1, \dotsc, f_p$ are polynomials. 
\end{definition}
Here we introduce our primary examples of polynomial maps and set systems. 

\begin{definition}\label{def:Veronese}
The \emph{degree $m$ Veronese map of $\Rl^d$} is the map $V_m^d: \Rl^d \to \Rl^{\binom{d+m}{m}-1}$ which maps $(x_1, \dotsc , x_d)$ to the vector (in some order) of all non-constant monomials of degree at most $m$ in the $d$ variables $x_1, \dotsc , x_d$.
\end{definition}

\begin{definition}\label{def:PH-system}
The \emph{degree $m$ homogeneous Veronese map of $\Rl^d$} is $HV_{m}^d: \Rl^d \to \Rl^{\binom{d+m-1}{m}}$ which maps $(x_1, \dotsc , x_d)$ to the vector (in some order) of all monomials of degree $m$ in the $d$ variables $x_1, \dotsc , x_d$.
\end{definition}

We will use the notation $(\Rl^d,\mathcal{P}_m^d)$ for the set system induced by the degree $m$ Veronese map of $\Rl^d$ and $(\Rl^d, \mathcal{H}_m^d)$ for the set system induced by the degree $m$ homogeneous Veronese map of $\Rl^d$. 

As a concrete example, consider the degree 2 Veronese map of $\Rl^2$. This is the map $V_2^2:\Rl^2 \to \Rl^5$ where $V_2^2(x,y) = (x^2,xy,y^2,x,y)$. The set system induced by this map is $(\Rl^2, \mathcal{P}_2^2)$. Its subsets consist of all regions of the plane determined by some conic section.

\subsection{On \(k\)-sets and \(k\)-facets for set systems induced by maps}


The natural notion of a $k$-set of a set system is a range that contains exactly $k$ points.
It is often more convenient to work with $k$-facets but it is not clear how to define them for set systems whose ranges lack a well-defined boundary. 
This motivates our restriction to set systems induced by maps, as their ranges have a well-defined boundary and interior.

\begin{definition}
Given a set system $(X,\F_\varphi)$ induced by a map $\varphi: X \to \Rl^p$ and a finite set $S \subset X$, an \emph{$\F_\varphi$-$k$-set} of $S$ is
a subset $A \subset S$ of size $k$ such that $\varphi(A)$ can be  strictly separated from $\varphi(S\setminus A)$ by a hyperplane. 
\end{definition}

\begin{definition}\label{def:kfacet}
Given a set system $(X,\F_\varphi)$ induced by a map $\varphi: X \to \Rl^p$ and a finite set $S \subset X$ such that $\varphi(S)$ is in general linear position, an \emph{$\F_\varphi$-$k$-facet} of $S$ is a subset $P$ of $p$ points from $S$, along with some orientation of the hyperplane $\aff \varphi(P)$, such that the subset $\varphi(P)$ along with the chosen orientation of $\aff \varphi(P)$ is a $k$-facet of $\varphi(S)$. 
\end{definition}


Observe that counting $\mathcal{F}_\varphi$-$k$-sets/facets simply amounts to counting $k$-sets/facets of point sets of the form $\varphi(S)$. Therefore, upper bounds for $e_k^{(d)}(n)$ and $a_k^{(d)}(n)$ immediately imply  non-trivial upper bounds on the number of $\F_\varphi$-$k$-sets/facets that a set of points may have:

\begin{proposition}\label{prop:embeddingbound}
Given a set system $(X,\F_\varphi)$ induced by a map $\varphi: X \to \Rl^p$, and a finite subset $S \subset X$ such that $\varphi(S)$ is in general linear position, the number of $\F_\varphi$-k-facets of $S$ is at most $e_k^{(p)}(n)$. 
\end{proposition}

For $\F_{\varphi}$-$k$-sets, we do not need to assume that $\varphi(S)$ is in general linear position since $k$-sets are defined for any point set whether or not it is in general linear position. 

\begin{proposition}\label{prop:embeddingbound2}
Given a set system $(X,\F_\varphi)$ induced by a map $\varphi: X \to \Rl^p$, the number of $\F_\varphi$-k-sets that a set of $n$ points in $X$ may have is at most $a_k^{(p)}(n)$. 
\end{proposition}
\makeproof{}{prop:embeddingbound2}{

In the case when $\varphi$ is not injective, $\varphi(S)$ may need to be considered as a multiset.
Therefore, we start with the observation that $a_k^{(p)}(n)$ is the maximum number of $k$-sets even for point sets which have repeated points, i.e., multisets of points. To see this, observe that perturbing a set of points can only increase the number of $k$-sets \cite{WagnerSurvey}. Therefore, if we start out with a multiset, we can perturb it slightly to create a set (in general linear position) with the same number of points and at least as many $k$-sets. Now, the number of $\F_\varphi$-$k$-sets of $S$ is equal to the number of $k$-sets of $\varphi(S)$, which is at most $a_k^{(p)}(n)$.
}

\ifsocg
\else
\subsection{Neighborly polytopes}\label{sec:neighborly}
For a set of $n$ points in convex position in the plane, the number of $k$-facets is precisely $n$ for all values of $k$. 
In $\Rl^3$, a similar result is true: the number of $k$-facets for a set $S$ of $n$ points in general position which form the vertex set of a 3-polytope is $2(k+1)n-4 \binom{k+2}{2}$, see \cite{WagnerSurvey}. 
There is no such result in dimension $d \ge 4$, i.e., convex position does not force a point set in $\Rl^d$ ($d \ge 4$) to have a specific number of $k$-facets. 
In fact, the $k$-set/$k$-facet problem for point sets in convex position in $\Rl^4$ is only slightly better understood than the problem for arbitrary point sets, see \cref{prop:convex}. 

However, if we assume that our point set is not only in convex position but is also \emph{neighborly}, then $e_k^{(d)}(S)$ is determined precisely by $\lvert S \rvert$. 
See \cite{WagnerSurvey,Ziegler} for an introduction to neighborly polytopes.
\begin{definition}
A polytope is $k$-\emph{neighborly} if any set of $k$ or fewer vertices forms a face. 
A $d$-polytope is \emph{neighborly} if it is $\lfloor d/2 \rfloor $-neighborly. 
\end{definition}

If we are talking about point sets instead of polytopes, we will say that a set of points is \emph{$k$-neighborly} (respectively, \emph{neighborly}) if it is the vertex set of a $k$-neighborly (respectively, neighborly) polytope. 

One way of producing neighborly point sets is by choosing a finite set of points on the \emph{moment curve} $M:= \{(x,x^2,\dotsc, x^d): x \in \Rl\}\subset \Rl^d$. 

\begin{proposition}[\cite{Clarkson87,DBLP:journals/dcg/ClarksonS89},{\cite[Proposition 4.1]{WagnerSurvey}}]\label{prop:Neighborly}
Let $S$ be a neighborly set of $n$ points in general linear position in $\Rl^d$. Then
\[
e_k(S) = 
\begin{cases*}
 2 \binom{k+\lceil d/2\rceil -1}{\lceil d/2 \rceil -1}\binom{n-k-\lceil d/2\rceil }{\lceil d/2 \rceil -1} 
 &if $d$ is odd \\[.1in]
 \binom{k+d/2-1}{d/2-1}\binom{n-k-d/2}{d/2}+\binom{k+d/2}{d/2}\binom{n-k-d/2-1}{d/2-1}
 &if $d$ is even.
\end{cases*}
\]
\end{proposition}

\fi
\section{Counting \texorpdfstring{$k$}{k}-facets via maps}\label{sec:count}
In this section we count $\mathcal{F}_{\varphi}$-$k$-facets when $\mathcal{F}_{\varphi}$ is a set system induced by a map. When the map $\varphi$ has certain properties, we can say more about the number of $\mathcal{F}_{\varphi}$-$k$-facets.

\subsection{Counting \(k\)-facets exactly}
It turns out that the maps associated to several families of polynomials we have discussed have the surprising property that they map generic point sets into the set of vertices of a neighborly polytope. 
Given such a map $\varphi$, we are able to exactly count the number of $\F_\varphi$-$k$-facets for point sets in general position. 

Before stating the new results, we recall a result of \cite{1676031,Ardila04} which served as motivation. 
A \emph{halving circle} of a point set of size $2n+1$ is a circle which has $3$ points on its boundary and $n-1$ points on either side. In the following theorem general position means that no three points are collinear and no four are concyclic. 
\begin{theorem}[\cite{1676031,Ardila04}]\label{thm:circles}
Any set $S$ of $2n + 1$ points in general position in the plane has exactly $n^2$ halving circles. More generally, for any $0 \le k \le 2n-2$, the number of circles that have 3 points of $S$ on their boundary and $k$ points on one side is exactly $2(k+1)(2n-k-1)$\footnote{This formula counts each halving circle twice, once for each orientation.}.
\end{theorem}
The proof of \cref{thm:circles} in \cite{Ardila04} is by a continuous motion argument. However, as noted there, it is possible to give a shorter proof using the method of maps as follows. The set system of all circles in the plane can be described as $(\Rl^2, \F_C)$ where $C:\Rl^2 \to \Rl^3$ is the map $C(x,y) = (x,y,x^2+y^2)$. Since $C$ maps generic point sets into convex position (on the surface of a paraboloid), \cref{thm:circles} follows from an application of the formula in \cref{prop:Neighborly} since any $3$-polytope is neighborly.

Now we define halving polynomials for other families of polynomials. Informally, for a finite set $S \subseteq \Rl^2$, a \emph{halving conic section} of $S$ is a conic section inequality having $5$ points on its boundary and half of the remaining points of $S$ in its interior. Unlike \cref{thm:circles} on the circle problem above, we count halving conic sections twice, once for each orientation. 
This is to be consistent with the standard definition of $k$-facets (\cref{def:kfacet}).
More precisely,
\begin{definition}
For a set $S$ of $2n+5$ points in $\Rl^2$, a \emph{halving conic section} is an $\F_V$-$n$-facet of $S$ where $V:=V_2^2$ is the degree 2 Veronese map of $\Rl^2$.
\end{definition}

\begin{definition}
For a set $S$ of $2n+m+1$ points in $\Rl^2$, a \emph{halving homogeneous polynomial of degree $m$} of $S$  is an $\mathcal{H}_m^2$-$n$-facet of $S$.  
\end{definition}
The halving case is a particular case of the more general problem of counting $k$-facets. It is generally believed that the maximum number of $k$-facets is maximized by the halving case, so it is considered the most important. However, we state our results counting $\mathcal{H}_m^2$-$k$-facets and $\F_V$-$k$-facets for arbitrary values of $k$.


\begin{theorem}\label{prop:ConicNeighborly}
Assume a finite set of points $S\subseteq \RR^2$ is in general linear position. Then the image of $S$ by the degree 2 Veronese map $V_2^2$ is neighborly. 
\end{theorem}
\makeproof{}{}{
First we verify that $V_2^2(S)$ is the set of vertices of a polytope. There is a bijection between conic sections passing through points of $S$ and hyperplanes passing through the images of those points by $V_2^2$. Therefore, for every point $v \in S$ we need to find a conic section inequality passing through $v$ and with all other points on one side. We can use an inequality of the form $(x-a)^2+ (y-b)^2 \le r$ and adjust the constants $a,b,r$ so that $(x-a)^2+ (y-b)^2 \le r$ defines a circle that contains $v$ on its boundary and has radius small enough so that no other point in $S$ is in the circle.  Now we show that $V_2^2(S)$ is neighborly. For every 2 points $v_1,v_2$ of $S$ we need to find a conic section inequality passing through those points and with all other points on one side. One way to accomplish this is to use the line $ax+by=c$ through $v_1$ and $v_2$. Then $(ax+by-c)^2\le 0$ is the required conic section inequality. 
}
In terms of the terminology defined in Section \ref{sec:generally}, the above result says that $V_2^2$ is a ``generally neighborly embedding''. The same result holds for the even degree homogeneous Veronese map of the plane. 
\begin{theorem}\label{prop:HomNeighborly}
Assume $m$ is even and $S\subseteq \RR^2$ is in general position, meaning that no two points of $S$ lie on a common line through the origin. Then the image of $S$ by $HV_{m}^2$ is neighborly. 
\end{theorem}
\makeproof{}{}{
The proof is similar to that of \cref{prop:ConicNeighborly}. For any set $\{v_1, \dotsc , v_{k}\}$ of $k \le m/2$ points of $S$ we need to find a degree $m$ homogeneous polynomial inequality which passes through all the $v_i$ and has all other points of $S$ on one side. Let $v_{k+1} , \dotsc, v_{m/2}$ be points in the plane that belong to no line passing through the origin and a point of $S$. For each $i$, $1 \le i \le m/2$, let $a_i x+b_i y=0$ be the line through the origin and $v_i$. Then $\prod_{i=1}^{m/2}(a_i x+b_i y)^2\le 0$ is a polynomial inequality with the required properties. 
}
The next results require us to strengthen our general position assumptions from \cref{prop:ConicNeighborly,prop:HomNeighborly}. 
\begin{definition}\label{def:genconic}
A set $S \subset \Rl^2$ is in \emph{general position with respect to conics} if $S$ is in general linear position and $V_2^2(S)$ is in general linear position. 
\end{definition}

\begin{definition}\label{def:genhom}
A set $S \subset \Rl^2$ is in \emph{general position with respect to degree $m$ homogeneous polynomials} if no two points in $S$ lie on a common line through the origin and $HV_{m}^2(S)$ is in general linear position. 
\end{definition}

\begin{theorem}\label{thm:coniccount}
Any set of $n$ points of $\Rl^2$ in general position with respect to conics has exactly $2 \binom{k+2}{2}\binom{n-k-3 }{2} $ $\F_V$-$k$-facets where $V:=V_2^2$ is the degree 2 Veronese map of the plane.

\end{theorem}
\makeproof{}{}{
Let $S \subset \Rl^2$ be a set of $n$ points in general position with respect to conics. 
There is a bijection between conic sections passing through 5 points of $S$ and hyperplanes passing through 5 points of $V(S)$. 
Furthermore, there is a bijection between $\F_V$-$k$-facets of $S$ and $k$-facets of $V(S)$. By \cref{prop:ConicNeighborly}, $V(S)$ is neighborly. 
Also, since $S$ is in general position with respect to conics, $V(S)$ is in general linear position. Therefore the number of $\F_V$-$k$-facets of $S$ is given by the formula from \cref{prop:Neighborly}.
}

\begin{theorem}\label{thm:homogeneouscount}
Assume $m$ is even. Any set of $2n+m+1$ points of $\Rl^2$ in general position with respect to degree $m$ homogeneous polynomials has exactly $2 \binom{k+m/2}{m/2}\binom{n-k-m/2-1 }{m/2} $ $\mathcal{H}_m^2$-$k$-facets. 
\end{theorem}
\makeproof{}{}{
Let $S \subset \Rl^2$ be a set of $2n+m+1$ points in general position with respect to degree $m$ homogeneous polynomials. 
As in the last proof, there is a bijection between $\mathcal{H}_m^2$-$k$-facets of $S$ and $k$-facets of $HV_{m}^2(S)$. By \cref{prop:HomNeighborly}, $HV_{m}^2(S)$ is neighborly. 
Also, since $S$ is in general position with respect to degree $m$ homogeneous polynomials, $HV_{m}^2(S)$ is in general linear position. 
Since $HV_m^2: \Rl^2 \to \Rl^{m+1}$, the formula from \cref{prop:Neighborly} in the case $d=m+1$ completes the proof.
}

\ifsocg
We can generate more set systems with the property that the number of $\F$-$k$-facets is independent of the configuration of points (as long as the points do not violate some general position assumption) by a lifting of the moment curve\lnote{preliminaries/appendix},
see \cref{sec:momentcurve}.
\else
\subsection{Lifting the moment curve}
We say that a set system $(X,\mathcal{F})$ has the exact count property if the number of $\mathcal{F}$-$k$-sets for any set of $n$ points in general position depends only on $n$ and $k$ and not on the configuration of points. We can generate many more set systems with the exact counting property by a lifting of the moment curve.

Let $f:\RR^d \to \RR$ be a function such that $\{(x,x') \in \Rl^{2d} : f(x)\neq f(x')\}$ is open and dense and let $g: \RR^d \to \RR$ be any function. 
We will say that a set $S = \{s_i\}_{i \in [n]}$ of $n$ points in $\Rl^d$ is in general position if $\prod_{i,j \in [n], i \neq j } \big(f(s_i)-f(s_j)\big) \neq 0$.\details{$f$ is injective on $S$. $f(s_i)$s are all distinct}
Note that this is a reasonable definition of general position since if we are considering $n$-point sets, then point sets (in $\Rl^{nd}$) in general position are open and dense in $\Rl^{nd}$. 
Assume that $m \geq 2$ is even. The map $\varphi: \Rl^d \to \Rl^{m+1}$ given by 
\begin{equation}\label{eq:lift}
\varphi(x) = \Bigl(f(x) , \bigl(f(x)\bigr)^2, \dotsc , \bigl(f(x)\bigr)^m, g(x)\Bigr)
\end{equation}
satisfies that, for any set $S$ of $n$ points in general position in $\Rl^d$, $\varphi(S)$ is neighborly.\footnote{This is saying that $\varphi$ is generally neighborly, using the terminology of \cref{sec:generally}. However, note that $\varphi$ may not be an embedding of $\Rl^d$. Furthermore, if $d>2$, any map $\varphi$ constructed as in \cref{eq:lift} cannot be an embedding of $\Rl^d$.}
To see this, note that, since $m$ is even, $\lfloor \frac{m+1}{2}\rfloor =m/2$. 
The projection of $\varphi(S)$ to the first $m$ coordinates is $m/2$-neighborly since it is a set of $n$ distinct points on the moment curve in $\Rl^m$. We claim that this implies that $\varphi(S)$ is $m/2$-neighborly as well: 
Let $\pi (\varphi(S))$ denote the projection to the first $m$ coordinates. 
The points of $\varphi(S)$ all project to distinct vertices of $\pi (\varphi(S))$. 
By neighborliness of $\pi (\varphi(S))$, every subset $F$ of at most $m/2$ vertices of $\pi(\varphi(S))$ forms a face. 
For any such face, there is a supporting hyperplane $H$. 
The preimage $\pi^{-1}(H)$ of $H$ under the projection $\pi$ is a hyperplane with normal having last coordinate 0. 
Moreover, $\pi^{-1}(H)$ is the supporting hyperplane for a face of $\varphi(S)$ formed by the lifted vertices $\pi^{-1}(F)$. 
This shows that $\varphi(S)$ is neighborly.

Given any admissible choice of functions $f,g$, the map above induces a set system $(\Rl^d, \F_\varphi)$ with the exact count property.


\fi

\subsection{Improved bounds for \(\mathcal{P}_m^d\)-\(k\)-facets and \( \mathcal{H}_m^d\)-\(k\)-facets}%
\label{sec:convex}

Results like \cref{thm:coniccount,thm:homogeneouscount} are not possible for any of the other polynomial set systems we have discussed. However, some progress can be made.

Recall that in \cref{prop:embeddingbound} we proved a non-trivial upper bound for the number of $\mathcal{F}$-$k$-facets where $(X,\mathcal{F})$ is any set system induced by a map. In this section we show how to improve this result for the set system $(\Rl^d,\mathcal{P}_m^d)$ for all values of $m$ and $d$ and for the set system $(\Rl^d,\mathcal{H}_m^d)$ for $m$ even.

We show that the maps which induce $(\Rl^d,\mathcal{P}_m^d)$ and $(\Rl^d,\mathcal{H}_{2m}^d)$, although not neighborly, still map into convex position with a high degree of neighborliness. Since these maps come up often in many fields, the following results may be useful in other contexts.  

\begin{theorem}\label{prop:veroneseneighborly}
For a finite set $S$ of points in $\Rl^d$, $V_m^d(S)$ is the set of vertices of an   $\ell$-polytope where $\ell  \le \binom{m+d}{m}-1$. If $V_{m/2}^d(S)$ is in general linear position and $m\geq 2$ is even, then $V_m^d(S)$ is a $\left(\binom{m/2+d}{m/2}-1\right)$-neighborly $\ell$-polytope.
\end{theorem}
\makeproof{}{prop:veroneseneighborly}{
For $v \in S$, choose coefficients $a_i$, $R$ so that $\{x \in \RR^m \suchthat \sum_{i=1}^m (x_i-a_i)^2 \le R \}$ is a ball with $v$ on its boundary and with radius small enough so that no points of $S \setminus v$ are inside. By the Veronese map $V^d_m$, this ball corresponds to a hyperplane in $\Rl^{\binom{m+d}{m}-1}$ containing $v$ and with all other points of $S$ on one side. This shows that $V_m^d(v)$ is a vertex of $\conv(V_m^d(S))$. 
For the second claim, let $T \subset S$, $|T| \le \binom{m/2+d}{m/2}-1$. Let $p(x)=1$ be a degree $m/2$ polynomial passing through each point of $T$ and no points of $S \setminus T$. 
To show that such a polynomial exists, recall we are assuming that $V_{m/2}^d(S)$ is in general position. 
Therefore, for any $|T|$ points in $V_{m/2}^d(S)$ there is a hyperplane passing through precisely those $|T|$ points. 
This hyperplane corresponds to a degree $m/2$ polynomial passing through each point of $T$ and no points of $S \setminus T$. 
Then $\bigl(p(x)-1\bigr)^2= 0$ is a polynomial surface which corresponds to a hyperplane in $\Rl^{\binom{m+d}{m}-1}$ which supports $\conv(T)$ as a face of $\conv(S)$.
}

\begin{theorem}\label{prop:homogeneousneighborly}
 Assume $m\geq 2$ is even. For a finite set $S$ of points in $\Rl^d$, $HV_{m}^d(S)$ is the set of vertices of an $\ell$-polytope where $\ell  \le \binom{m+d-1}{m}$. 
 If $HV_{m/2}^d(S)$ is in general position, meaning no hyperplane through the origin in the image space of $HV_{m/2}^d$ contains more than $\left(\binom{m/2+d-1}{m/2}-1\right)$ points of $HV_{m/2}^d(S)$, then $HV_{m}^d(S)$ is a $\left(\binom{m/2+d-1}{m/2}-1\right)$-neighborly $\ell$-polytope.
\end{theorem}
\makeproof{}{prop:homogeneousneighborly}{
For $v \in S$, let $H = \{x \in \Rl^d : a\cdot x = 0\}$ be a plane through the origin which contains $v$ and contains no other point of $S$. Then $(a\cdot x)^m \le 0$ is a degree $m$ homogeneous polynomial inequality which, by the homogeneous Veronese map $HV^d_m$, corresponds to a hyperplane in $\Rl^{\binom{m+d-1}{m}}$ containing $v$ and with all other points of $S$ on one side. This shows that $HV_{m}^d(v)$ is a vertex of $\conv(HV_{m}^d(S))$. For the second claim, let $T \subset S$, $|T| \le \left(\binom{m/2+d-1}{m/2}-1\right)$. Let $p(x)=0$ be a degree $m/2$ homogeneous polynomial passing through each point of $T$ and no points of $S \setminus T$. To show that such a polynomial exists, recall we are assuming that $HV_{m/2}^d(S)$ is in general position. 
Therefore, for any $|T|$ points in $HV_{m/2}^d(S)$ there is a hyperplane passing through the origin and precisely those $|T|$ points. And this hyperplane corresponds to a degree $m/2$ homogeneous polynomial passing through each point of $T$ and no points of $S \setminus T$. Then $(p(x))^2 = 0$ is a degree $m$ homogeneous polynomial surface which corresponds to a hyperplane in $\Rl^{\binom{m+d-1}{m}}$ which supports $\conv(T)$ as a face of $\conv(S)$. 
}

These $k$-neighborliness results are of interest to us because convex position is a special case of the $k$-set/facet problem for which we can improve the best known upper bound:
\begin{theorem}\label{prop:convex}
For a set $S$ of $n$ points in convex position in $\Rl^d$, $e_k(S) \le (n/d)e_k^{(d-1)}(n-1)$. 
\end{theorem}
\makeproof{}{}{
Let $v \in S$. Choose a hyperplane $H$ containing $v$ and with all other points of $S$ on one side of it. Choose another hyperplane $H'$ parallel to $H$ and with all points of $S$ between $H$ and $H'$. 
Let $S'$ be the stereographic projection (using $v$ as the ``pole'') of $S \setminus v$ onto $H'$. 
We claim that the number of $k$-facets of $S$ containing $v$ is equal to the number of $k$-facets of $S'$ (as a subset of $H'$, a $(d-1)$-dimensional subspace). 

Assume that $\conv(v_1 ,\dotsc , v_{d-1},v)$ is a $k$-facet of $S$. We claim that for each $s \in S$, the stereographic projection $s'$ is on the positive side of $\aff(v_1', \dotsc ,v_{d-1}')$ if and only if $s$ is on the positive side of $\aff(v_1, \dotsc , v_{d-1},v)$. This is seen to be true by observing that $\aff(s,v)$ does not intersect $\aff(v_1, \dotsc , v_{d-1},v)$ anywhere other than the point $v$. This shows that $\conv(v_1', \dotsc , v_{d-1}')$ is a $k$-facet of $S'$. For the converse, assume that $\conv(v_1', \dotsc ,v_{d-1}')$ is a $k$-facet of $S'$. Then $\conv(v_1, \dotsc ,v_{d-1},v)$ is a $k$-facet of $S$ for the same reason as above. 

Since $S'$ lies in a hyperplane, it can have at most $e_k^{(d-1)}(n-1)$ $k$-facets. 
Performing this projection on each point of $S$ and noticing that every $k$-facet is counted $d$ times shows the desired result.
}
As far as we know, the best known bound for $k$-facets of $n$-point sets in $\Rl^4$ in convex position is the same as for general point sets, which is $O(n^2k^{2-2/45})$ \cite{4DBOUND}. 
In $\Rl^3$ the best known bound for general point sets is $O(nk^{2-1/2})$ \cite{3DBOUND,3DBOUNDJ}. 
This combined with \cref{prop:convex} gives a bound of $O(n^2k^{2-1/2})$ for the number of $k$-facets of point sets in convex position in $\Rl^4$.

An argument similar to the proof of \cref{prop:convex} shows the following generalization (we state it without proof):
\begin{proposition}\label{prop:nbound}
For a set $S$ of $n$ points in $m$-neighborly position in $\Rl^d$, $e_k(S) \le \frac{\binom{n}{m}e_k^{(d-m)}(n-m)}{\binom{d}{m}}$.
\end{proposition}

In dimensions higher than 4, the best known bound for $k$-facets is $e_{k}^{(d)}(n) = O(n^{d-\epsilon_d})$ where $\epsilon_d = (4d-3)^{-d}$ \cite{HIGHDBOUND}. Because of the fast decay of the constant $\epsilon_d$, \cref{prop:nbound} gives an improvement in the best known upper bound which depends on the degree of neighborliness of the point set in question.

\cref{prop:nbound} can be used to improve the bounds for $\mathcal{P}_m^d$-$k$-facets and $ \mathcal{H}_{2m}^d$-$k$-facets as follows. Recall that the set system $(\Rl^d, \mathcal{P}_m^d)$ is induced by the map $V_m^d$ and $(\Rl^d, \mathcal{H}_{2m}^d)$ is induced by the map $HV_{2m}^d$. Therefore \cref{prop:veroneseneighborly,prop:homogeneousneighborly} along with \cref{prop:nbound} give an improvement in the bound.


\section{Limits of the neighborliness argument}\label{sec:limits}
In this section we study under what conditions exact count phenomenon can occur for arbitrary set systems induced by a map. 

The crucial observation that allows us to exactly count $\mathcal{F}$-$k$-facets for the conic sections and even degree homogeneous polynomials on the plane is that the maps which induce these set systems map generic point sets to neighborly point sets. 
In this section we define a \emph{generally neighborly embedding} to be an embedding that maps generic point sets to neighborly point sets (\cref{def:generallyneighborly}). 
The moment curve map is an example of a generally neighborly embedding of $\Rl^1$. Moving up one dimension, the degree 2 Veronese map of the plane $V_2^2$ shows that a generally neighborly embedding of the plane exists (by \cref{prop:ConicNeighborly}). 
We conjecture that generally neighborly embeddings of $\Rl^d$ do not exist for $d>2$.  In order to provide support for this conjecture, in \cref{sec:variety} we prove a closely related result about algebraic varieties. More evidence that the conjecture may be true is provided in \cref{sec:evidence}.

Apart from its relation to the $\mathcal{F}$-$k$-facets, the problem of determining the existence of generally neighborly embeddings is interesting in its own right. Our definition of generally neighborly embeddings is similar to and inspired by Micha Perles' definition of \emph{neighborly embeddings} which we discuss in \cref{subsec:perles}.

\subsection{Generally neighborly embeddings}\label{sec:generally}


\begin{definition}[embedding]
An \emph{embedding} is a map which is a homeomorphism onto its image.
\end{definition}

\begin{definition}[generally $k$-neighborly embedding]\label{def:generallyneighborly}
Let $\varphi: \Rl^d \to \Rl^p$ be an embedding. 
For each $n \in \Nl$, let $G_n \subset \Rl^{dn}$ consist of all configurations of $n$ points in $\Rl^d$ which are mapped to $k$-neighborly sets by $\varphi$. Then $\varphi$ is \emph{generally $k$-neighborly} if $G_n$ contains a set that is open and dense in $\Rl^{dn}$ for all $n$. A generally $\lfloor \frac{p}{2}\rfloor$-neighborly embedding is called a \emph{generally neighborly embedding}.
\end{definition}
We choose ``open and dense'' in \cref{def:generallyneighborly} for concreteness and readability. For part of our discussion (in particular, \cref{prob:dim} below), it may be reasonable to substitute it by an alternative version of a generic property as discussed in \cref{sec:generic}.

Observe that the even degree homogeneous Veronese map of the plane, i.e. $HV_m^2$, is not an embedding because it is not injective. Since the homogeneous Veronese map is one of our prime examples throughout we need to justify why we are now only talking about embeddings. The reason is that all of the polynomial maps we have considered have the property that they are an embedding of some open subset of Euclidean space. 
Thus, for the purposes of this section it suffices to assume that our maps are embeddings.

The main question concerning generally $k$-neighborly embeddings is:

\begin{problem}\label{prob:dim}
What is the smallest dimension $p:= p_g(k,d)$ of the image space for which a generally
$k$-neighborly embedding $\varphi: \Rl^d \to \Rl^p$ exists? 
\end{problem}
%
%

\begin{theorem}\label{thm:embedding}
There exists a generally $k$-neighborly embedding of $\Rl^d$ into $\Rl^{2k+d-1}$ and so $p_g(k,n)  \le  2k+d-1$.
\end{theorem}
\makeproof{}{}{
Consider the embedding $\varphi:\Rl^d \to \Rl^{2k+d-1}$ defined by
\[
\varphi (x_1,x_2, \dotsc, x_d) = (x_1,x_1^2, x_1^3, \dotsc, x_1^{2k},x_2, \dotsc , x_d).
\]
For each $n \in \Nl$, let $G_n \subset \Rl^{dn}$ consist of all configurations of $n$ points in $\Rl^d$ such that no two points in the configuration have the same $x_1$-coordinate. One can verify that $G_n$ is open and dense in $\Rl^{dn}$. 
For any $n$, let $S \in G_n$ be some configuration of $n$ points in $G_n$. 

To show that $\varphi(S)$ is $k$-neighborly, let $v_1, \dotsc , v_k$ be $k$ points from $S$. 
Consider in the domain of the embedding, $\RR^d$, the surface
\begin{equation}\label{equ:surface}
\prod_{i=1}^k (x_1-v_{i1})^2 = 0.
\end{equation}
By expanding \cref{equ:surface} we see that this surface corresponds via $\varphi$ to a hyperplane $H$ in $\Rl^{2k+d-1}$.
Note that $v_1, \dotsc , v_k$ satisfy \cref{equ:surface} and all other points in $S$ satisfy $\prod_{i=1}^k (x_1-v_{i1})^2 > 0$.
Using $\varphi$, we get that $H$ is a face-defining hyperplane that makes $\varphi(v_1), \dotsc , \varphi(v_k)$ a face of $\conv(\varphi(S))$.
Therefore, $\varphi(S)$ is $k$-neighborly. 
This shows that $p_g(k,d) \le 2k+d-1$.
}
We believe that the bound in the above theorem is actually tight.

\begin{conjecture}\label{Conjecture:conjecture}
$p_g(k,d) = 2k+d-1$. 
\end{conjecture}

Observe that for $d \ge 3$, if $\varphi: \Rl^d \to \Rl^p$ is a generally $k$-neighborly embedding then, according to \cref{Conjecture:conjecture}, $p \ge 2k+2$. This means the conjecture implies that generally neighborly embeddings of $\Rl^d$ do not exist for $d \ge 3$. 

In the context of the $k$-set problem, this would mean that set systems like the conic sections do not exist in dimension $d \ge 3$. More precisely, it would imply that, for $d \ge 3$, there is no set system $(\Rl^d, \mathcal{F})$ induced by an embedding (of $\Rl^d$) for which $\F$-$k$-facets can be counted by using \cref{prop:Neighborly}.


\subsection{Generally neighborly manifolds} 

Here we define generally $k$-neighborly manifolds which are, for our purposes, equivalent to generally $k$-neighborly embeddings. The sense in which they are equivalent is made precise below.

\begin{definition}\label{def:generallyneighborlymanifold}
A manifold $M \subset \Rl^p$ is \emph{generally $k$-neighborly} if the set $G_n \subset M^n$ of configurations of $n$ points on $M$ which are $k$-neighborly contains a set that is open and dense in $M^n$ for all $n$.  A generally $\lfloor \frac{p}{2}\rfloor$-neighborly manifold is called a \emph{generally neighborly manifold}.
\end{definition}

We ask the same question for manifolds as we did for embeddings:

\begin{problem}\label{prob:dimM}
What is the smallest dimension $p$ of the ambient space in which a generally $k$-neighborly $d$-manifold $M \subset \Rl^p$ exists? 
\end{problem}

Observe that an open subset of a generally $k$-neighborly $d$-manifold is still generally $k$-neighborly. 
Therefore, in the context of \cref{prob:dimM} it suffices to assume that the manifold $M$ is (globally) homeomorphic to $\Rl^d$, that is $M = \varphi(\Rl^d)$ for some embedding $\varphi$. This observation, along with the following proposition, shows that \cref{prob:dimM} is equivalent to \cref{prob:dim}.  
\begin{proposition}
An embedding $\varphi: \RR^d \to \RR^p$ is generally $k$-neighborly if and only if $M:=\varphi(\RR^d)$ is a generally $k$-neighborly $d$-manifold. 
\end{proposition}
\makeproof{}{}{
If $\varphi: \Rl^d \to \Rl^p$ is generally $k$-neighborly, then for each $n$ the set of configurations of $n$ points which are mapped by $\varphi$ to $k$-neighborly point sets contains a set $O$ that is open and dense in $\Rl^{dn}$. Let $N \subset M^n$ be the set of $k$-neighborly configurations of $n$-points in $M^n$. The set $N$ contains $(\varphi \times \dotsb \times \varphi)( O)$ which is open and dense since $\varphi \times \dotsb \times \varphi$ is a homeomorphism. Therefore $M : = \varphi(\Rl^d)$ is a generally $k$-neighborly $d$-manifold. The other direction is similar.
}

\subsection{On \(k\)-neighborly algebraic varieties}\label{sec:variety}

Recall that \cref{prob:dimM} asks for the minimal dimension of the ambient space in which a generally $k$-neighborly $d$-manifold exists. In this section we introduce an algebraic version of this problem. We ask a similar question with manifolds replaced by algebraic varieties. 
The problem we study concerns \emph{generally $k$-neighborly  algebraic varieties}. Before defining such varieties, we review some necessary background material from algebraic geometry. 

\subsubsection{Preliminaries on classical algebraic geometry}\label{sec:alggeo}

Given polynomials $f_1, f_2, \dotsc , f_r \in \Rl[x_1, \dotsc, x_p]$, the \emph{affine algebraic variety} defined by the $f_i$'s is the set 
$$
V(\Cx) = \{x \in \Cx^p : f_1(x) = 0, \dotsc, f_r(x)=0\}. 
$$

For a given variety $V(\Cx) \subset \Cx^p$, we are mainly interested in the subset $V(\Rl):=V(\Cx) \cap \Rl^p$ of real points. We refer to $V(\Rl)$ as an \emph{affine real algebraic variety}. Unless otherwise stated, by an \emph{algebraic variety} or simply \emph{variety} we mean an affine real algebraic variety.

Since we are mainly interested in the set of real points, we will write $V$ for $V(\Rl)$ and will write $V(\Cx)$ to indicate when the complex points are also considered. See \cite{Bochnak,Harris} for more definitions from real algebraic geometry.

\begin{definition}
A variety $V \subset \Rl^p$ (resp. $V(\Cx) \subset \Cx^p)$ is \emph{non-degenerate} if it is not contained in any hyperplane in $\Rl^p$ (resp. $\Cx^p$). 
\end{definition}

\begin{definition}
An algebraic variety is \emph{irreducible} if it cannot be written as the union of two proper algebraic subvarieties. 
\end{definition}

A well known fact from algebraic geometry is that every variety can be written as a finite union of irreducible components. 

We will use facts about the smooth and singular points of a variety. Suppose $V(\Cx)$ is a variety and the ideal of $V(\Cx)$ is generated by the polynomials $f_1, \dotsc, f_r$. The \emph{smooth} points of $V(\Cx)$ are those points where the Jacobian matrix of the $f_i$'s has maximal rank.
A \emph{singular point} of a variety is a point that is not smooth.

\begin{definition}
We use $V_{\text{sm}}$ to denote the set of smooth points of an algebraic variety $V$, $V_{\text{sing}}$ is the set of singular points.
\end{definition}

Given a real variety $V \subset \Rl^p$, let $V_{\Cx}$ denote the smallest complex variety which contains $V$. It is well known that $V_{\Cx}$ is unique and furthermore that there is a bijection between irreducible components of $V$ and irreducible components of $V_{\Cx}$, see \cite{WhitneyElemReal}. Note that $V(\Cx)$ is not always equal to $V_{\Cx}$. However, it is a useful fact that whenever $V$ contains a smooth real point, $V$ is Zariski dense in $V(\Cx)$ and so $V(\Cx) = V_{\Cx}$, see \cite[Section 2.8]{Bochnak}. We will always assume that our varieties contain smooth real points. 

By the \emph{dimension} of a real algebraic variety $V$ we will mean the dimension of $V(\Cx)$. There is another notion of dimension for real algebraic varieties. 

\begin{definition}
The \emph{real dimension} of a real algebraic variety $V$ is the maximal integer $d$ such that there is a homeomorphism of $[0,1]^d$ into some subset of $V$.
\end{definition}

The real dimension of $V$ does not always equal the dimension of $V$. However, the real dimension of $V$ is never more than the dimension of $V$ (see \cite[Proposition 2.8.14]{Bochnak}) and if $V$ contains a smooth real point, then these dimensions do agree. This is because around any smooth real point of a $d$-dimensional variety $V \subset \Rl^p$, there is a neighborhood which is a smooth $d$-dimensional submanifold of $\Rl^p$ \cite[Proposition 3.3.11]{Bochnak}.

In one of the proofs in \cref{sec:weakvar} we will also consider projective varieties, see \cite{Harris} for background on projective algebraic geometry. We use $\mathbb{P}^p(K)$ to denote $p$-dimensional projective space over the field $K = \Cx$ or $K = \Rl$.

\subsubsection{Generally $k$-neighborly algebraic varieties}
\begin{definition}\label{def:genvariety}
Let $V \subset \Rl^p$ be an irreducible real algebraic variety with a smooth real point. For each $n$, let $G_n \subset V_{\text{sm}}^n$ consist of all configurations of $n$ points on $V_{\text{sm}}$ which are $k$-neighborly. Then $V$ is \emph{generally $k$-neighborly} if $G_n$ contains a set that is open and dense in $V_{\text{sm}}^n$ for all $n$. It is \emph{generally neighborly} if it is generally $\lfloor p/2 \rfloor$-neighborly. 
\end{definition}

We make one clarifying remark regarding the above definition. One could replace $V_{\text{sm}}$ everywhere in the above definition with $V$. However, only requiring the property to hold for the smooth points strengthens our results below and does not change the proofs. Another reason for only considering the smooth points is the following. Loosely speaking, a generally $k$-neighborly algebraic variety $V$ is supposed to be a variety such that every generic configuration of points on $V$ is $k$-neighborly. A generic configuration of points should never contain non-smooth points, so the non-smooth points should be ignored when defining generally $k$-neighborly algebraic varieties. 

The question we are dealing with in this section is the following. 
\begin{problem}\label{probvariety}
What is the smallest dimension $p:=p_{g,V}(k,d)$ of the ambient space in which a generally
$k$-neighborly $d$-dimensional algebraic variety $V \subset \Rl^p$ exists? 
\end{problem}

Observe that the image of the map $\varphi$ in \cref{thm:embedding} is a $d$-dimensional generally $k$-neighborly variety in $\Rl^{2k+d-1}$. This shows that $p_{g,V}(k,d) \le 2k+d-1$.

We will prove the following result which completely resolves \cref{probvariety}.

\begin{theorem}\label{prop:lowerbound}
 Let $V \subset \Rl^p$ be a generally $k$-neighborly $d$-dimensional algebraic variety\footnote{Note that ``generally $k$-neighborly'' requires $V$ to be irreducible and contain a smooth real point.}. 
Then $p \geq 2k+d-1$.
\end{theorem}
\cref{prop:lowerbound} combined with \cref{thm:embedding} show that $p_{g,V}(k,d) = 2k+d-1$. 
In order to prove \cref{prop:lowerbound} we will first establish a connection between generally $k$-neighborly varieties and \emph{weakly $k$-neighborly} sets. 
Weak neighborliness is a more usable property that holds for all subsets of points, not just those satisfying a general position assumption. This connection is established in \cref{sec:weaklysets}. 
The proof of \cref{prop:lowerbound} is then completed in \cref{sec:weakvar}. 

Below we list some more examples of generally $k$-neighborly algebraic varieties. 
\begin{example}
The image of the degree 2 Veronese map of the plane is a generally 2-neighborly 2-dimensional algebraic variety in $\Rl^5$. 
\end{example}
\begin{example}
The image of the map $\varphi$ from the proof of \cref{thm:embedding} is a generally $k$-neighborly $d$-dimensional algebraic variety in $\Rl^{2k+d-1}$. 
\end{example}
\begin{example}
The \emph{moment curve} is a generally neighborly 1-dimensional algebraic variety. 
The same is true of any \emph{order $d$ curve} which is also an algebraic variety. 
\end{example}

\begin{example}
\cref{prop:lowerbound} shows that generally neighborly $d$-dimensional algebraic varieties do not exist for $d \ge 3$. 
\end{example}

\subsubsection{Weakly \(k\)-neighborly sets}\label{sec:weaklysets}
It turns out that all generally $k$-neighborly algebraic varieties satisfy a weaker neighborliness property that holds for all subsets of points (not just those satisfying a general position assumption). 
We call this property \emph{weakly $k$-neighborly} (\cref{def:weakly}). 
In the proof of \cref{prop:lowerbound}, we only need to use the fact that $V_{\text{sm}}$ is weakly $k$-neighborly. 
In this section we prove some lemmas concerning weakly $k$-neighborly sets and the relationship between generally $k$-neighborly manifolds/varieties and weakly $k$-neighborly sets. 

\begin{definition}\label{def:weakly}
A set $S \subseteq \RR^p$ is \emph{weakly $k$-neighborly} if for any set $T$ of $k$ points from $S$, there exists a closed halfspace $H$ with boundary $\bd(H)$ such that $S \subset H$ and $T \subset \bd(H)$. 
\end{definition}

We will now show that a generally $k$-neighborly algebraic variety or manifold is weakly $k$-neighborly. 
In order to do so, we first show that every finite subset of such a variety or manifold is weakly $k$-neighborly. 
We actually state and prove a stronger result (\cref{lem:generallyweakly}) that only uses as assumption the ``dense'' part of ``open and dense'' in the definition of generally $k$-neighborly.
We then establish a compactness property of weakly $k$-neighborly sets. 
The property is that an arbitrary subset of $\RR^p$ is weakly $k$-neighborly if and only if every finite subset is weakly $k$-neighborly. 
These two results together establish that generally $k$-neighborly algebraic varieties and manifolds are weakly $k$-neighborly.

\begin{lemma}\label{lem:generallyweakly}
%
Let $M \subset \Rl^p$ be a manifold or the set of smooth points of an algebraic variety.
If the set $N \subset M^n$ of configurations of $n$ points on $M$ which are $k$-neighborly is dense in $M^n$ for all $n$, then every finite set of points on $M$ is weakly $k$-neighborly. 
\end{lemma}
\makeproof{}{lem:generallyweakly}{
Assume not, so that there exists some finite set $S \subset M$  and a set $T$ of $k$ points from $S$ such that no closed halfspace contains $S$ and contains $T$ on its boundary. 
This means that $\aff(T) \cap \relint \conv(S \setminus T) \neq \emptyset$ (from the separating hyperplane theorem \cite[Theorem 1.3.8]{MR3155183}). 
We can pick $|S|$ small open balls in $\Rl^p$ as follows. For each $t \in T$, let $B_t$ be a ball centered at $t$ and for each $s \in S \setminus T$, let $A_s$ be a ball centered at $s$. The radii of the balls can be chosen small enough so that any collection of points consisting of one point from each $B_t$ and one point from each $A_s$ has the property that the affine hull of the points from the $B_t$ intersects the relative interior of the convex hull of the points from the $A_s$. This means that any such configuration of points is not $k$-neighborly. Therefore, the subset of $\Rl^{p|S|}$ of configurations of points of size $|S|$ on $M$ which are not $k$-neighborly contains $\prod_{t \in T } (B_t \cap M) \times \prod_{s \in S \setminus T } (A_s \cap  M)$ which is an open subset of $M^{|S|}$. 
Therefore, the set of configurations of $|S|$ points on $M$ which are $k$-neighborly is not dense in $M^{|S|}$.
}


\begin{proposition}[compactness]\label{prop:compactness}
Let $S \subseteq \RR^p$ be a (possibly infinite) set.
Then for any $k \geq 1$ we have that $S$ is weakly $k$-neighborly if and only if every finite subset of $S$ is weakly $k$-neighborly.
\end{proposition}
\makeproof{}{prop:compactness}{
Fix $k \geq 1$.
The ``only if'' direction is clear.
We will now prove the ``if'' direction.
Let $T \subseteq S$ be a set of $k$ points.
Let $\mathcal{U} = \{ U \subseteq S \suchthat U \supseteq T \text{ and $U$ is finite}\}$.
For $U \subseteq S$ such that $U \supseteq T$, we will define $N(U) \subseteq S^{p-1}$ to be the set of unit outer normals to possible halfspaces $H$ such that $T \subseteq \bd(H)$ and $U\subseteq H$.
More precisely, let 
\[
N(U) = \{a \in S^{p-1} \suchthat 
(\forall x,y \in T) a \cdot x = a \cdot y \text{ and } (\forall x \in T) (\forall y \in U) a \cdot x \geq a \cdot y \}.
%
\]
Clearly $N(U)$ is closed.
Let $\mathcal{V} \subseteq \mathcal{U}$ be any finite subfamily.
Then $\cap_{U \in \mathcal{V}} N(U) = N(\cup_{U \in \mathcal{V}} U) \neq \emptyset$ by assumption.
We have established that $\{N(U)\}_{U \in \mathcal{U}}$ is a family of closed sets with the finite intersection property in compact space $S^{p-1}$.
This implies $\cap_{U \in \mathcal{U}} N(U) \neq \emptyset$.
We also have $N(S) = N(\cup_{U \in \mathcal{U}} U) = \cap_{U \in \mathcal{U}} N(U) \neq \emptyset$.
That is, there is a halfspace $H$ such that $T \subseteq \bd(H)$ and $S\subseteq H$.
As $T \subseteq S$ was arbitrary, this completes the proof.
}

We require two more lemmas concerning intersections and weak separation of convex sets in $\Rl^p$.

Two sets $Q,R \subseteq \RR^p$ can be \emph{weakly separated} if there exist a non-zero $a\in \RR^p$ and $t \in \RR$ such that $Q\subseteq \{x \in \RR^p \suchthat a \cdot x \leq t \}$ and $R\subseteq \{x \in \RR^p \suchthat a \cdot x \geq t \}$. We say that the hyperplane $\{x \in \RR^p \suchthat a \cdot x =t\}$ \emph{weakly separates} $Q$ from $R$. 
This separation is said to be \emph{proper} if $Q$ and $R$ are not both contained in $\{x \in \RR^p \suchthat a \cdot x = t \}$.
\begin{lemma}[Radon-type theorem]\label{lem:radon}
Let $P$ be a set of $p+2$ points in $\RR^p$ in general linear position. 
Then there is a partition $Q,R$ of $P$ into two non-empty sets so that $\relint \conv Q \cap \relint \conv R \neq \emptyset$.
\end{lemma}
\makeproof{}{lem:radon}{
Let $P = \{ q_1, \dotsc, q_{p+2}\}$. 
$P$ is affinely dependent and therefore there exist $\lambda_1, \dotsc, \lambda_{p+2}$ such that $\sum_{i=1}^{p+2} \lambda_i q_i = 0$, $\sum_{i=1}^{p+2} \lambda_i = 0$ and at least one $\lambda_i$ is non-zero.
Because of the general position assumption, all $\lambda_i$ are non-zero.
Let $I = \{ i \suchthat \lambda_i >0 \}$, $J = \{ i \suchthat \lambda_i < 0\}$.
Both $I$ and $J$ are non-empty.
By dividing $\lambda_i$s by $\sum_{i \in I} \lambda_i$ we can assume without loss of generality that $\sum_{i \in I} \lambda_i = -\sum_{i \in J} \lambda_i = 1$.
Let $Q = \{q_i \suchthat i \in I \}$, $R = \{q_i \suchthat i \in J \}$.
Let $q = \sum_{i \in I} \lambda_i q_i \in \relint \conv Q$, $r = -\sum_{i \in J} \lambda_i q_i \in \relint \conv R$.
We have $q=r$, which completes the proof.
}

\begin{lemma}\label{lem:separation}
Let $Q, R \subseteq \RR^p$ be disjoint sets. 
Suppose 
$\aff (Q \cup R) = \RR^p$
and 
\(
\relint \conv Q \cap \relint \conv R \neq \emptyset
\).
Then $Q, R$ cannot be weakly separated.
\end{lemma}
\makeproof{}{lem:separation}{
Assume $Q,R$ can be weakly separated. 
If the separation is not proper then $\aff(Q \cup R) \neq \Rl^p$. 
If the separation is proper, then by the separating hyperplane theorem \cite[Theorem 1.3.8]{MR3155183}, $\relint \conv Q \cap \relint \conv R = \emptyset$.
}

\subsubsection{Weakly \(k\)-neighborly varieties}\label{sec:weakvar}

In the previous section we established the connection between generally $k$-neighborly varieties/manifolds and weakly $k$-neighborly sets. In this section, we use this connection to prove \cref{prop:lowerbound}. Given an (real) algebraic variety $V$ of dimension $d$ such that $V_{\text{sm}}$ is weakly $k$-neighborly, we prove a sharp lower bound on the dimension of the ambient space.
\begin{theorem}\label{thm:lowerboundweak}
Assume $V \subset \Rl^p$ is a non-degenerate $d$-dimensional irreducible real algebraic variety with a smooth real point. If $V\setminus U$ is weakly $k$-neighborly for some proper closed subvariety $U$, then $p \ge 2k+d-1$.
\end{theorem}

Before proving \cref{thm:lowerboundweak} we prove the special case of algebraic curves and then generalize to higher dimensional varieties.

\begin{lemma}\label{lem:lowerboundcurve}
Assume $C \subset \Rl^p$ is a non-degenerate irreducible real algebraic curve with a smooth real point. If $C \setminus U$ is weakly $k$-neighborly for some proper closed subvariety $U$, then $p \ge 2k$. 
\end{lemma}
\begin{proof}
First we will show that one can find arbitrarily large point sets in general linear position on $C \setminus U$. In order to accomplish this, first observe that $C \setminus U$ is non-degenerate. Indeed, if it were the case that $C \setminus U$ is contained in a hyperplane $H$, then we would have $C = (C \cap H) \cup U$ which is impossible since $C$ is irreducible. 

Now assume that $S$ is a set of $j$ points in general linear position on $C \setminus U$. If it were not possible to find another point $s$ such that $S \cup  \{s\}$ is in general linear position, it would have to be the case that $C \setminus U$ is contained in the union of all hyperplanes spanned by $p$ points from the set $S$. Let $\mathcal{H}$ be the collection of all hyperplanes spanned by $p$ points in $S$. Note that $\mathcal{H}$ is finite. We have
\[
C = \bigg(\bigcup_{H \in \mathcal{H}} \bigl((C \setminus U) \cap H\bigr) \bigg) \cup U.
\]
Since $C \setminus U$ is non-degenerate the above formula would be a representation of $C$ as the union of proper subvarieties. This is impossible since $C$ is irreducible. Therefore we know that $C \setminus U$ is not contained in the union of all hyperplanes spanned by points in $S$ and so we can always find $s$ so that $S \cup \{s\}$ is in general linear position. It follows that we can find arbitrarily large point sets in general linear position on $C \setminus U$.

Let $P$ be a set of $p+2$ points in general linear position on $C \setminus U$. By \cref{lem:radon,lem:separation}, there is a partition of $P$ into non-empty sets $Q$ and $R$ so that $Q$ and $R$ cannot be weakly separated.

However, because $C \setminus U$ is weakly $k$-neighborly, we know that for any set $T$ of $k$ points on $C\setminus U$ there exists a closed halfspace which contains $C\setminus U$ and contains $T$ in its boundary. 
In other words, any such $T$ can be weakly separated from $C\setminus U$. Therefore, it must be that $k< \min(|Q|,|R|)$, i.e. $k \le \min(|Q|,|R|)-1$. 
Now since $\min(|Q|,|R|)  \le \lfloor\frac{p+2}{2}\rfloor$, we have that $k \le \lfloor\frac{p+2}{2}\rfloor - 1$ and so $p \ge 2k$.
\end{proof}

The idea of the proof for higher dimensional varieties is to take successive hyperplane sections in order to reduce to the case of curves. So we first need to establish the following lemma concerning hyperplane sections of varieties.

\begin{lemma}\label{lem:hypsection}
Assume that $V \subset \Rl^p$ is a non-degenerate $d$-dimensional ($d \ge 2$) irreducible variety with a smooth real point and that $U$ is some proper closed subvariety of $V$. 
Then for any given open ball $B$ in $\conv( V\setminus U)$ there exists a hyperplane $H$ such that $H \cap B \neq \emptyset$ and $V \cap H$ is a non-degenerate, irreducible, $(d-1)$-dimensional variety with a smooth real point which is not contained in $U$. 
\end{lemma}
\begin{proof}
We identify the set of hyperplanes in $\Rl^p$ with a proper subset of $\mathbb{P}^p(\Rl)$.
This identification works as follows. 
We identify a hyperplane $a_0 +a_1x_1 + \dotsb + a_px_p=0$
in $\Rl^p$ with the point with homogeneous coordinates $(a_0,a_1, \dotsc , a_p)$ in $\mathbb{P}^p(\Rl)$. 
Therefore, the set of hyperplanes in $\Rl^p$ is identified with $\mathbb{P}^p(\Rl) \setminus \{(a_0, a_1, a_2, \dotsc , a_p) \suchthat a_1 = a_2 = \cdots =a_p =0\}$.

Recall from \cref{sec:alggeo} that around any smooth real point of $V$ there is a neighborhood which is a smooth $d$-dimensional submanifold of $\Rl^p$. 
Let $T \subset \mathbb{P}^p(\Rl)$ be the set of hyperplanes $H$ for which there exists a smooth point of $V \setminus U$ and a neighborhood of that smooth point which has nonempty transversal intersection with $H$. 
We know that $T$ is open and that for any $H \in T$, $H \cap V$ contains an open
subset which is a smooth $(d-1)$-dimensional submanifold of $\Rl^p$ \cite[Sections 1.5 and 1.6]{Pollack}. 
Let $H$ be any hyperplane in $T$. 
Clearly the dimension of $H \cap V$ is at most $d-1$.
We claim that $H \cap V$ is a variety of dimension precisely $d-1$ and that it contains a smooth real point.
Recall from \cref{sec:alggeo} that if a variety has real dimension $d$ then it has dimension at least $d$.
Therefore, the dimension of $H \cap V$ is $d-1$. 
Now we show that $H \cap V$ contains a smooth real point. 
Indeed assume not, so that $H \cap V(\Cx)$ contains no smooth real points.%
\footnote{By a minor abuse of notation, $H \cap V(\Cx)$ is the complex variety defined by the polynomials defining $V$ along with the polynomial defining $H$. 
Similarly, $H \cap V$ is the real variety defined by the polynomials defining $V$ along with the polynomial defining $H$. 
So $H \cap V$ is the real part of $H \cap V(\Cx)$.} 
This means that $H \cap V$ is contained in the set of singular points of the $(d-1)$-dimensional variety $H \cap V(\Cx)$. 
By \cite[Proposition 3.3.14]{Bochnak}, the set of singular points of $H \cap V(\Cx)$ is a variety of dimension at most $d-2$. 
Since $H \cap V$ has real dimension $d-1$ and is contained in a variety of dimension at most $d-2$ this is a contradiction. 
So $H \cap V$ contains smooth real points. 

Let $I$ be the subset of $T$ consisting of hyperplanes $H$ such that $H \cap V$ is irreducible and non-degenerate. 
We will show that $I$ is open and dense in the standard topology in $T$. Let $\overline{V(\Cx)} \subset \mathbb{P}^p(\Cx)$ be the projective closure of $V(\Cx)$. 
This means that $V(\Cx) = \overline{V(\Cx)}\cap \{x_0 \neq 0\}$. 
Because $V$ contains a smooth real point, $V$ is Zariski dense in $V(\Cx)$, that is, $V(\Cx)$ is the smallest complex variety containing $V$, see \cite[Section 2.8]{Bochnak}\details{alternatively, Theorem 5.1 in \url{https://arxiv.org/pdf/1606.03127.pdf}}. Therefore, since $V$ is irreducible, by \cite[Lemma 7]{WhitneyElemReal}, $V(\Cx)$ is irreducible. 
Because $V(\Cx)$ is irreducible, it is then a standard fact that $\overline{V(\Cx)}$ is irreducible.
By projective duality, we can identify the set of all projective hyperplanes with real coefficients in $\mathbb{P}^p(\Cx)$ with $\mathbb{P}^p(\Rl)$. 
Let $I' \subset \mathbb{P}^p(\Rl)$ consist of all projective hyperplanes with real coefficients that have irreducible and non-degenerate intersection with $ \overline{V(\Cx)}$. 
By \cite[Theorem 18.10]{Harris}, $I'$ is Zariski open and dense.%
\footnote{Theorem 18.10 in \cite{Harris} says that the set of hyperplanes which intersect $\overline{V(\Cx)}$ in a non-degenerate irreducible variety is the complement of a proper subvariety of $\mathbb{P}^p(\Cx)$. The intersection of a proper subvariety of $\mathbb{P}^p(\Cx)$ with $\mathbb{P}^p(\Rl)$ is a proper subvariety of $\mathbb{P}^p(\Rl)$. So $I'$ is Zariski open dense in $\mathbb{P}^p(\Rl)$.}

We will now show that the fact that $I'$ is Zariski open and dense implies that $I$ is also Zariski open and dense in $T$.  
Given a hyperplane $H$ defined by $a_0 + a_1x_1 + \dotsb + a_px_p=0$ in $\Rl^p$, there is a corresponding projective hyperplane $\overline{H}$ defined by $a_0x_0 + a_1x_1 + \dotsb + a_px_p=0$ which is the homogenization of $H$. 
Let $H$ be a hyperplane in $T$ and assume that the homogenization $\overline{H}$ is in  $I'$. 
We claim that this implies that $H \in I$. 
To establish this claim, we need to verify that $H \cap V$ is irreducible and nondegenerate. 
Observe that $H \cap V(\Cx) = \overline{H} \cap \overline{V(\Cx)} \cap \{x_0  \neq 0\}$ i.e., $H \cap V(\Cx)$ is an open subset of $ \overline{H} \cap \overline{V(\Cx)}$. 
It is a standard fact that a nonempty open subset of an irreducible space is irreducible and dense.
Therefore $H \cap V(\Cx)$ is irreducible. Now, since $H \cap V(\Cx)$ is Zariski dense in $\overline{H} \cap \overline{V(\Cx)}$, if $H \cap V(\Cx)$ were contained in a hyperplane, $\overline{H} \cap \overline{V(\Cx)}$ would be as well.
So $H \cap V(\Cx)$ is nondegenerate. 
We have shown that $H \cap V(\Cx)$ is non-degenerate and irreducible.
We claim that because the section $H \cap V(\Cx)$ contains smooth real points, then the non-degeneracy and irreducibility of $H \cap V(\Cx)$ implies non-degeneracy and irreducibility of $H \cap V$.
This relies on the fact mentioned above that if an irreducible variety $W$ contains a smooth real point, then the set of real points is Zariski dense in $W(\Cx)$.
This means that $H \cap V$ is Zariski dense in $H \cap V(\Cx)$. Now by \cite[Lemma 7]{WhitneyElemReal}, irreducibility of $H \cap V(\Cx)$ implies that $H \cap V$ is irreducible. 
Finally, observe that because $H \cap V$ is Zariski dense in $H \cap V(\Cx)$, if $H \cap V$ is degenerate, $H \cap V(\Cx)$ must be as well. 
Thus $H \cap V$ is non-degenerate.

Recall that $T$ is open (in the standard metric topology) and that $I'$ is Zariski open and dense in $\mathbb{P}^p(\Rl)$ which implies that $I'$ is open and dense in the standard topology. We established that $I$ contains $I' \cap T$.  Therefore, we have shown that $I$ is open and dense in the standard topology in $T$.

Let $O \subset \mathbb{P}^p(\Rl)$ be the set of hyperplanes which intersect $B$. The set $O$ is open (in the standard topology). 

We claim that all of this means that $I \cap T \cap O$ is non-empty. 
To establish this, we need to verify that $T \cap O$ is non-empty. 
To find $H \in T \cap O$, let $v$ be a smooth point of $V \setminus U$ and let $F$ be a $(p-2)$-flat having non-empty transversal intersection with a neighborhood of $v$. 
Such a flat exists because $V$ has real dimension at least 2. 
For any point $p$ in $B$, the hyperplane $\aff(F \cup p)$ is in $T \cap O$.

Therefore, $T \cap O$ is open and non-empty.
Since $I$ is open and dense in $T$, it follows that $I \cap T \cap O \neq \emptyset$. 
We claim that any $H \in I \cap T \cap O$ completes the proof. To show this, it remains to show that for any $H \in I \cap T \cap O$, $H \cap V$ has a smooth real point which is not in $U$. We already know that $H \cap V$ has smooth real points and that $H \cap V$ is not contained in $U$. So assume for a contradiction that all the smooth real points are contained in $U$. Then letting $S$ denote the singular points of $H\cap V$, we have that $H \cap V = (H \cap V \cap U)\cup S$ is the decomposition of $H \cap V$ as the union of two proper closed subvarieties, a contradiction to irreducibility of $H\cap V$. 
\end{proof}

The lemma established above allows us to generalize \cref{lem:lowerboundcurve} to higher dimensional varieties.

\begin{proof}[Proof of \cref{thm:lowerboundweak}]
First we show that by making repeated applications of \cref{lem:hypsection}, we can inductively construct a $(p-d+1)$-flat $L$ that intersects $\interior \conv (V \setminus U)$ and such that $L \cap V$ is a non-degenerate irreducible algebraic curve with a smooth real point which is not contained in $U$. 
Say we have some flat $F$ that intersects $\interior \conv (V \setminus U)$ and such that $F \cap V$ is a non-degenerate irreducible $d'$-dimensional $(d' \ge 2)$ variety with a smooth real point not contained in $U$. 
Since $F \cap V$ is not contained in $U$, $F \cap U$ is a proper subvariety of $F \cap V$. 
Let $B$ be an open ball in $F$ that is contained in $\interior \conv (V \setminus U)$. 
By \cref{lem:hypsection}, there exists a hyperplane $H$ in $F$ that intersects $B$ and such that $H \cap V$ is a non-degenerate irreducible $(d'-1)$-dimensional variety with a smooth real point not contained in $F \cap U$ and hence not contained in $U$. 
Since $B$ is contained in $\interior \conv (V \setminus U)$, $H$ intersects $\interior \conv (V \setminus U)$. 
We can repeat this process until we obtain a $(p-d+1)$-flat $L$ that intersects $\interior \conv (V \setminus U)$ and such that $C:=L \cap V$ is a non-degenerate irreducible algebraic curve with a smooth real point not contained in $U$.

Now we will show that $C \setminus U$ is weakly $k$-neighborly in $L$, meaning\footnote{Note that a subset of $\RR^p$ whose affine hull is a proper subset of $\RR^p$ is automatically weakly $k$-neighborly according to \cref{def:weakly}, so the requirement here is that halfspace $H$ in that definition is a halfspace of $L$.} that $C \setminus U$ is weakly $k$-neighborly as a subset of its affine hull $L$. 

Because $V \setminus U$ is weakly $k$-neighborly, we know that for any set $T$ of $k$ points on $C\setminus U$ there exists a closed halfspace $H$ (in $\Rl^p$) which contains $C \setminus U$ and contains $T$ in $\bd(H)$. 
In other words, any such $T$ can be weakly separated from $C \setminus U$. 
Notice that we are talking about weak separation in $\Rl^p$, while we are really interested in weak separation in $L$. 
We claim that our assumption that $L$ intersects the interior of the convex hull of $V \setminus U$ allows us to pass from weakly separating hyperplanes in $\Rl^p$ to weakly separating hyperplanes in $L$. 
Indeed, the fact that $L$ intersects $\interior \conv(V \setminus U)$ means that any closed halfspace $H$ satisfying $V \setminus U \subset H$  cannot contain $L$ in its boundary. 
Therefore $\bd(H) \cap L$ is a proper hyperplane in $L$. 
To summarize, given a set $T$ of $k$ points on $C\setminus U$, by the assumption that any finite set on $V \setminus U$ is weakly $k$-neighborly, there exists a hyperplane $\bd(H)$ that weakly separates $T$ from $V \setminus U$. 
This and the way we chose $L$ allows us to conclude that $\bd(H) \cap L$ is a hyperplane in $L$ that weakly separates $T$ from $C\setminus U$. Therefore, $C \setminus U$ is weakly $k$-neighborly in $L$. Since $C$ is not contained in $U$, $C  \cap U$ is a proper subvariety. By \cref{lem:lowerboundcurve}, $p-d+1 \ge 2k$.
\end{proof}

We can now prove the lower bound for generally $k$-neighborly varieties.
\ifdcg
\begin{proof}[Proof of \cref{prop:lowerbound}]
\else
\begin{proof}[Proof of \cref{prop:lowerbound}]
\fi
We can without loss of generality assume that $V$ is non-degenerate since otherwise we could consider $V$ in $\aff(V)$. By \cref{lem:generallyweakly}, every finite set of points on $V_{\text{sm}}$ is weakly $k$-neighborly. Therefore by \cref{prop:compactness}, $V_{\text{sm}}$ is weakly $k$-neighborly. Since $V_{\text{sing}}$ is a proper closed subvariety (see \cite[Proposition 3.3.14]{Bochnak}), by \cref{thm:lowerboundweak}, $p \ge 2k+d-1$. 
\ifdcg
\end{proof}
\else
\end{proof}
\fi

\subsection{Additional evidence}\label{sec:evidence}
We give some additional comments on the validity of \cref{Conjecture:conjecture}. 
Although we could not resolve the conjecture, the previous result on algebraic varieties is evidence that it is likely true. In the following we provide more evidence for the conjecture by showing that any manifold violating the conjecture would have to have a property that appears fairly restrictive to us. 
\begin{proposition}\label{prop:weird}
If a set $M \subset \Rl^p$ is weakly $k$-neighborly then for any set $S$ of $2k$ points in general linear position in $M$, $\aff(S) \cap M$ is contained in the union of all hyperplanes supporting facets of the simplex $\conv S$. 
\end{proposition}
\makeproof{}{}{
Assume not, so that $S$ is a set of $2k$ points in $M$ such that $\aff(S)\cap M$ is not contained in the union of all hyperplanes supporting facets of $\conv(S)$. 
Then identifying $\aff(S)$ with $\Rl^{2k-1}$, we can find a set $S'$ of $2k+1$ points in general position in $\Rl^{2k-1}$ which are weakly $k$-neighborly. 
However, by \cref{lem:radon,lem:separation}, there is a partition $Q,R$ of $S'$ into two non-empty sets so that $Q,R$ cannot be weakly separated. 
Since $\min(Q,R) \le \lfloor \frac{2k+1}{2} \rfloor= k$, this is a contradiction to weakly $k$-neighborliness of $S'$. 
}
\Cref{prop:weird} is inspired by the following illustrative example: 
the possibility of a non-degenerate 2-neighborly curve in $\RR^3$. 
By non-degeneracy, pick 4 affinely independent points on the curve. 
Then by \cref{prop:weird} the curve would have to be contained in the union of the 4 hyperplanes defined by any 3 of those points.

Assume $M \subset \Rl^p$ is a weakly $k$-neighborly embedded $d$-manifold and $S$ any set of $2k$ points in general position on $M$. Notice that if $p<2k+d-1$, we would expect $\aff(S)$ to intersect $M$ in a manifold of dimension $1$ or greater. However, the previous Proposition implies that $M \cap \aff S$ is contained in a finite number of hyperplanes which is not true of most embedded 1-manifolds. 

One approach to proving \cref{Conjecture:conjecture} appears to be showing that, in fact, there is no non-degenerate $d$-manifold in $\Rl^{2k+d-2}$ satisfying the conclusion of \cref{prop:weird}.

\ifsocg
\else
\subsection{Neighborly embeddings}\label{subsec:perles}

Our definition of generally $k$-neighborly embeddings is similar to the concept of \emph{$k$-neighborly embeddings} introduced by Perles in 1982 and studied by Kalai and Wigderson \cite{KalaiWigderson08}.

An embedding of a $d$-dimensional manifold $M$ into $\Rl^p$ is \emph{$k$-neighborly} if for every $k$ points on the embedding of $M$ there is a hyperplane $H$ that contains the $k$ points and such that all remaining points of the embedded manifold are on the (strictly) positive side of $H$. 

Requiring that an embedding be $k$-neighborly is clearly stronger than requiring that it be generally $k$-neighborly. That is, if an embedding is $k$-neighborly then it is generally $k$-neighborly. However, the reverse implication is certainly not true. For example, the degree two Veronese embedding $V_2^2(x,y) = (x,y,x^2,xy,y^2)$ is only a \emph{1}-neighborly embedding while it is a generally \emph{2}-neighborly embedding. 

In 1982 Perles posed the following problem concerning neighborly embeddings.

\begin{problem}
\normalfont What is the smallest dimension $p(k,d)$ of the ambient space in which a $k$-neighborly $d$-dimensional manifold exists?
\end{problem}

As in the case of generally neighborly embeddings, for the purposes of this question, it suffices to assume that $M = \Rl^d$. Kalai and Wigderson proved
\begin{theorem}[\cite{KalaiWigderson08}]\label{thm:Kalai}
$k(d+1) \le p(k,d) \le 2k(k-1)d$. 
\end{theorem}
Improving the bounds in \cref{thm:Kalai} appears to be difficult compared to the case of generally $k$-neighborly embeddings where we were able to conjecture a precise formula for $p_g(k,d)$. 

Comparing the two definitions, $k$-neighborly embeddings appear to be the more natural and fundamental class of embeddings to investigate. However, there may be some applications for which the notion of generally $k$-neighborly embeddings is more appropriate. For example, the authors of \cite{KalaiWigderson08} were interested in neighborly embeddings in part because they may lead us to important examples of $k$-neighborly polytopes. In particular, by picking points on the embedded manifold, one may produce  $k$-neighborly nonsimplicial polytopes perhaps with other interesting properties. Since both types of embeddings produce $k$-neighborly polytopes, our version may be more useful in this context as it is less restrictive.

\fi

\section{Conclusion and open questions}\label{sec:conclusion}
The main question left open concerns the validity of \cref{Conjecture:conjecture}. Of course resolving the conjecture in full generality would be ideal, but it may be worthwhile to focus on other special cases instead. For example, one might be able to prove the conjecture for analytic manifolds or for smooth manifolds using tangency properties. 


Another open question concerns \emph{lower bounds} for the algebraic $k$-set problem. In analogy with the original $k$-set problem we ask: 
Is there a polynomial map of the plane into some higher dimensional space which induces a natural set system $(\RR^2,\mathcal{F})$ for which the maximum number of $\mathcal{F}$-$k$-sets for a set of $n$ points is $n^\ell e^{\Omega (\sqrt{\log n})}$ and $o(n^{\ell+1})$ for some integer $\ell\ge 2$? Or just $\Omega(n^\ell\log n)$ and $o(n^{\ell+1})$?
A candidate that we do not fully understand in this context is the map $(x,y) \mapsto (x,y,xy)$ (or, equivalently up to a linear transformation, $(x,y) \mapsto (x,y,x^2-y^2)$).

\begin{acknowledgements}
We would like to thank Nina Amenta, Eric Babson, Imre B\'{a}r\'{a}ny, Kenneth Clarkson, Jes\'us De Loera, Tamal Dey, Igor Pak, Rom Pinchasi, Anne Schilling and Adam Sheffer for helpful discussions. We would also like to thank the anonymous reviewers whose comments improved the structure and presentation of the paper. 
This material is based upon work supported by the National Science Foundation under Grant No.\ DMS-1440140 while the second author was in residence at the Mathematical Sciences Research Institute in Berkeley, California, during the Fall 2017 semester. 
This material is also based upon work supported by the National Science Foundation under Grants CCF-2006994, CCF-1657939 and CCF-1422830.
\end{acknowledgements}

\ifsocg
\else
\bibliographystyle{abbrv}
\fi

\bibliography{b}

\begin{thebibliography}{10}

\bibitem{HIGHDBOUND}
N.~Alon, I.~B\'{a}r\'{a}ny, Z.~F\"{u}redi, and D.~J. Kleitman.
\newblock Point selections and weak {$\epsilon$}-nets for convex hulls.
\newblock {\em Combin. Probab. Comput.}, 1(3):189--200, 1992.

\bibitem{Ardila04}
F.~Ardila.
\newblock The number of halving circles.
\newblock {\em The American Mathematical Monthly}, 111(7):586--591, 2004.

\bibitem{Bochnak}
J.~Bochnak, M.~Coste, and M.-F. Roy.
\newblock {\em Real Algebraic Geometry}.
\newblock Springer, 1998.

\bibitem{DBLP:conf/compgeom/BuzagloHP08}
S.~Buzaglo, R.~Holzman, and R.~Pinchasi.
\newblock On \emph{s}-intersecting curves and related problems.
\newblock In M.~Teillaud, editor, {\em Proceedings of the 24th {ACM} Symposium
  on Computational Geometry, College Park, MD, USA, June 9-11, 2008}, pages
  79--84. {ACM}, 2008.

\bibitem{DBLP:conf/compgeom/ChevallierFSS14}
N.~Chevallier, A.~Fruchard, D.~Schmitt, and J.~Spehner.
\newblock Separation by convex pseudo-circles.
\newblock In S.~Cheng and O.~Devillers, editors, {\em 30th Annual Symposium on
  Computational Geometry, SOCG'14, Kyoto, Japan, June 08 - 11, 2014}, page 444.
  {ACM}, 2014.

\bibitem{Chevallier}
N.~Chevallier, A.~Fruchard, D.~Schmitt, and J.~Spehner.
\newblock Separation by convex pseudo-circles.
\newblock {\em Discrete {\&} Computational Geometry}, 2020.

\bibitem{Clarkson87}
K.~L. Clarkson.
\newblock New applications of random sampling in computational geometry.
\newblock {\em Discrete {\&} Computational Geometry}, 2:195--222, 1987.

\bibitem{DBLP:journals/dcg/ClarksonS89}
K.~L. Clarkson and P.~W. Shor.
\newblock Application of random sampling in computational geometry, {II}.
\newblock {\em Discrete {\&} Computational Geometry}, 4:387--421, 1989.

\bibitem{Dey1997ImprovedBO}
T.~K. Dey.
\newblock Improved bounds on planar k-sets and k-levels.
\newblock {\em Proceedings 38th Annual Symposium on Foundations of Computer
  Science}, pages 156--161, 1997.

\bibitem{DBLP:journals/dcg/Dey98}
T.~K. Dey.
\newblock Improved bounds for planar k-sets and related problems.
\newblock {\em Discrete {\&} Computational Geometry}, 19(3):373--382, 1998.

\bibitem{MR0363986}
P.~Erd\H{o}s, L.~Lov\'{a}sz, A.~Simmons, and E.~G. Straus.
\newblock Dissection graphs of planar point sets.
\newblock In {\em A survey of combinatorial theory ({P}roc. {I}nternat.
  {S}ympos., {C}olorado {S}tate {U}niv., {F}ort {C}ollins, {C}olo., 1971)},
  pages 139--149, 1973.

\bibitem{Pollack}
V.~Guillemin and A.~Pollack.
\newblock {\em Differential topology}.
\newblock AMS Chelsea Publishing, Providence, RI, 2010.
\newblock Reprint of the 1974 original.

\bibitem{Harris}
J.~Harris.
\newblock {\em Algebraic geometry}, volume 133 of {\em Graduate Texts in
  Mathematics}.
\newblock Springer-Verlag, New York, 1995.
\newblock A first course, Corrected reprint of the 1992 original.

\bibitem{KalaiWigderson08}
G.~Kalai and A.~Wigderson.
\newblock Neighborly embedded manifolds.
\newblock {\em Discrete {\&} Computational Geometry}, 40(3):319--324, 2008.

\bibitem{1676031}
D.-T. Lee.
\newblock On k-nearest neighbor \uppercase{V}oronoi diagrams in the plane.
\newblock {\em IEEE Transactions on Computers}, C-31(6):478--487, June 1982.

\bibitem{lovasz1971number}
L.~Lov{\'a}sz.
\newblock On the number of halving lines.
\newblock {\em Ann. Univ. Sci. Budapest, Etovos, Sect. Math.}, 14:107--108,
  1971.

\bibitem{4DBOUND}
J.~Matou\v{s}ek, M.~Sharir, S.~Smorodinsky, and U.~Wagner.
\newblock {$k$}-sets in four dimensions.
\newblock {\em Discrete Comput. Geom.}, 35(2):177--191, 2006.

\bibitem{MR2405686}
G.~Nivasch.
\newblock An improved, simple construction of many halving edges.
\newblock In {\em Surveys on discrete and computational geometry}, volume 453
  of {\em Contemp. Math.}, pages 299--305. Amer. Math. Soc., Providence, RI,
  2008.

\bibitem{MR3155183}
R.~Schneider.
\newblock {\em Convex bodies: the {B}runn-{M}inkowski theory}, volume 151 of
  {\em Encyclopedia of Mathematics and its Applications}.
\newblock Cambridge University Press, Cambridge, second expanded edition, 2014.

\bibitem{3DBOUND}
M.~Sharir, S.~Smorodinsky, and G.~Tardos.
\newblock An improved bound for \emph{k}-sets in three dimensions.
\newblock In S.~Cheng, O.~Cheong, P.~K. Agarwal, and S.~Fortune, editors, {\em
  Proceedings of the Sixteenth Annual Symposium on Computational Geometry,
  Clear Water Bay, Hong Kong, China, June 12-14, 2000}, pages 43--49. {ACM},
  2000.

\bibitem{3DBOUNDJ}
M.~Sharir, S.~Smorodinsky, and G.~Tardos.
\newblock An improved bound for \emph{k}-sets in three dimensions.
\newblock {\em Discret. Comput. Geom.}, 26(2):195--204, 2001.

\bibitem{Shawe-Taylor:2004:KMP:975545}
J.~Shawe-Taylor and N.~Cristianini.
\newblock {\em Kernel Methods for Pattern Analysis}.
\newblock Cambridge University Press, New York, NY, USA, 2004.

\bibitem{DBLP:conf/compgeom/Toth00}
G.~T{\'{o}}th.
\newblock Point sets with many \emph{k}-sets.
\newblock In S.~Cheng, O.~Cheong, P.~K. Agarwal, and S.~Fortune, editors, {\em
  Proceedings of the Sixteenth Annual Symposium on Computational Geometry,
  Clear Water Bay, Hong Kong, China, June 12-14, 2000}, pages 37--42. {ACM},
  2000.

\bibitem{Tothlower}
G.~T{\'{o}}th.
\newblock Point sets with many \emph{k}-sets.
\newblock {\em Discrete {\&} Computational Geometry}, 26(2):187--194, 2001.

\bibitem{WagnerSurvey}
U.~Wagner.
\newblock {$k$}-sets and {$k$}-facets.
\newblock In {\em Surveys on discrete and computational geometry}, volume 453
  of {\em Contemp. Math.}, pages 443--513. Amer. Math. Soc., Providence, RI,
  2008.

\bibitem{WhitneyElemReal}
H.~Whitney.
\newblock Elementary structure of real algebraic varieties.
\newblock {\em Ann. of Math.}, (2):545--556, 1957.

\bibitem{Ziegler}
G.~M. Ziegler.
\newblock {\em Lectures on polytopes}.
\newblock Springer-Verlag, New York, 1995.

\end{thebibliography}

\ifsocg
\appendix

\section{Lifting the moment curve}\label{sec:momentcurve}

\fi

\end{document}